\documentclass[11pt]{amsart}

\numberwithin{equation}{section}
\usepackage{mathrsfs,txfonts}
\usepackage[dvips]{graphicx}
\usepackage{latexsym,bm}
\usepackage{amsfonts}
\usepackage{indentfirst}
\usepackage{amsthm}
\usepackage{cases}
\usepackage{pifont}
\usepackage{amsmath}
\usepackage{amssymb}

\allowdisplaybreaks
\newtheorem{lem}{\quad\textbf{\Large Lemma}}[section]
\newtheorem{thm}[lem]{\quad\textbf{\Large Theorem}}

\def\squarebox#1{\hbox to #1{\hfill\vbox to #1{\vfill}}}

\usepackage{enumerate}
\begin{document}
\title[Semi-finite proof of triple product identity]
{A semi-finite proof of Jacobi's triple product identity }
\author{
Jun-Ming Zhu }
\address{Department of Mathematics, Luoyang Normal University,
Luoyang City, Henan Province 471022, China}
\email{junming\_zhu@163.com}
\thanks{Keywords: Jacobi's triple product identity;  Euler's $q$-exponential function; semi-finite.
\\\indent MSC (2010):   33D15, 11F27.
% \\\indent Supported by the
% National Science Foundation of China
% (Grant No. 11371184).
}

\begin{abstract}
Jacobi's triple product identity is proved from one of Euler's
$q$-exponential functions in an elementary way.
% As its application, an identity on $q$-gamma function is obtained.
\end{abstract}
 \maketitle
%%%%%%%%%%%%%%%%%%%%%%%%%%%%%%%%%%%%%%%%%%%%%%%%%%%%%%%%%%%%%%%%%%%%%%%%%%%%%

\section{Introduction}
We  suppose $|q|<1$ and, as usual, the $q$-Pochhammer symbols are
defined respectively by
\begin{equation*}\label{def}
(a;q)_\infty=\prod_{n=0}^{\infty}(1-aq^{n}) \qquad \mbox{and}\qquad
(a;q)_n=\frac{(a;q)_\infty}{(aq^n;q)_\infty}
\end{equation*}
for any integer $n$. Jacobi's triple product identity is one of the
most important series-product identities. We state it in the
following theorem.
\begin{thm}\label{tri} For any complex number $z\neq0$, we have
\begin{equation}\label{jacobi}
(q;q)_\infty(-q/z;q)_\infty(-z;q)_\infty=\sum_{n=-\infty}^{\infty}q^{n(n-1)/2}z^n.
\end{equation}
\end{thm}
Bellman said in \cite[p. 42]{cf} that there are no simple proofs
known of the complete result. After that, Andrews \cite{and} proved
Theorem \ref{tri} using both of Euler's $q$-exponential functions,
which are, respectively,
\begin{equation}\label{E1}
(-z;q)_\infty=\sum_{n=0}^{\infty}\frac{q^{n(n-1)/2}z^n}{(q;q)_n},~~
z\in \mathbf{C},
\end{equation}
and
\begin{equation*}\label{E2}
\frac{1}{(-z;q)_\infty}=\sum_{n=0}^{\infty}\frac{(-1)^n
z^n}{(q;q)_n},~~ |z|<1.
\end{equation*}
In this short note, we will prove Jacobi's triple product identity
only using \eqref{E1}. This seems also a proof of the complete
triple product identity required by Bellman. We use a method  which
is often called  semi-finite method (see \cite{cf}).

\section{The proof}
\begin{proof} For a nonnegative integer $m$ and $z\neq0$,
\begin{eqnarray}
\sum_{n=-m}^{\infty}\frac{q^{n(n-1)/2}z^n}{(q^{m+1};q)_n}&=&\sum_{n=0}^{\infty}\frac{q^{(n-m)(n-m-1)/2}z^{n-m}}{(q^{m+1};q)_{n-m}}\label{E1x}\\
&=&\frac{q^{m(m+1)/2}z^{-m}}{(q^{m+1};q)_{-m}}\sum_{n=0}^{\infty}\frac{q^{n(n-1)/2}{(zq^{-m})}^n}{(q;q)_n}\notag\\
&=&q^{m(m+1)/2}z^{-m}(q;q)_m(-zq^{-m};q)_\infty\notag\qquad\mbox{(by \eqref{E1})}\\
&=&q^{m(m+1)/2}z^{-m}(q;q)_m(-zq^{-m};q)_m(-z;q)_\infty\notag\notag\\
&=&(q;q)_m(-q/z;q)_m(-z;q)_\infty.\notag
\end{eqnarray}
Note that   $\frac{1}{(q^{m+1};q)_{-n}}=(q^{m+1-n};q)_n=0$ when
$n>m$.
 The  left hand side of \eqref{E1x} can be rewritten as
\begin{eqnarray*}\label{E1xx}
%\sum_{n=0}^{\infty}\frac{q^{n(n-1)/2}z^n}{(q^{m+1};q)_n}+\sum_{n=1}^{\infty}\frac{q^{n(n+1)/2}z^{-n}}{(q^{m+1};q)_{-n}}
%=
\sum_{n=0}^{\infty}\frac{q^{n(n-1)/2}z^n}{(q^{m+1};q)_n}+\sum_{n=1}^{\infty}(q^{m+1-n};q)_nq^{n(n+1)/2}z^{-n}.
\end{eqnarray*}
 Restricting $z$ in any compact subset
of $0<|z|<\infty$  and letting $m\rightarrow\infty$ in \eqref{E1x},
we get \eqref{jacobi}. By analytic continuation, the restriction on
$z$ may be relaxed.
\end{proof}

{\bf Acknowledgment.} The author wishes to thank the editor and the
referee for their valuable comments and advices. The author is
partially supported by the
 National Science Foundation of China
(Grant No. 11371184).

%\begin{acknowledgment}{Acknowledgment.}
% The authors wish to thank the editor and the referee for
%their valuable comments and advices. The author is partially
%supported by the
% National Science Foundation of China
%(Grant No. 11371184).
%\end{acknowledgment}

%

\end{document}